\documentclass[twoside,11pt]{amsart} 
\usepackage{enumitem}
\usepackage{epsfig}
\usepackage{latexsym} 
\usepackage{amsmath} \usepackage{amssymb}

 \keywords{ primes, gaps, prime constellations, Eratosthenes sieve}

\subjclass{11N05, 11A41, 11A07}

\newtheorem{theorem}{Theorem}[section]
\newtheorem{lemma}[theorem]{Lemma}

\newdimen\epsfxsize
\newdimen\epsfysize

\newcommand {\biggap}     {\makebox[0.2 in]{}}   
   
\newcommand{\lil}   {\scriptstyle }

\newcommand {\pml}[1]  {{#1}^{\#}}

\newcommand{\pgap}   {{\mathcal G}}

\begin{document}

\title{Addendum:  models for gaps $g=2p_{1}$}

\date{25 Sept 2023}

\author{Fred B. Holt}
\address{fbholt62@gmail.com; https://www.primegaps.info}

\begin{abstract}
This is an addendum to previous work on the models of populations and relative populations of small gaps $g$ across
stages of Eratosthenes sieve.  If we have the initial conditions in the cycle of gaps $\pgap(\pml{p_0})$, we can exhibit
exact population models for all gaps $g < 2 p_1$.  For gaps beyond this threshold $2 p_1$, we could not be certain of the
counts of the driving terms for a gap $g$ of various lengths.

Here we extend this work by introducing the exact population models for $g=2p_1$.  We are able to get this one additional 
case in a general form.  Using the initial conditions from $\pgap(\pml{p_0})$ we advance the model one time to obtain
exact counts for the driving terms for $g=2p_1$ in $\pgap(\pml{p_1})$.  This iteration uses a different system matrix than the
usual one.  After this special iteration we can apply the usual dynamics to obtain the exact population model for the 
gap $g=2p_1$ across all stages of Eratosthenes sieve.

\end{abstract}

\maketitle

\section{Setting}

This is an addendum to previous work \cite{FBHSFU, FBHPatterns} on exact models for the populations and relative populations
of gaps $g$ across stages of Eratosthenes sieve.  At each stage of Eratosthenes sieve, there is a cycle of gaps $\pgap(\pml{p})$
of length $\phi(\pml{p})$ (number of gaps in the cycle) and span $\pml{p}$ (sum of the gaps in the cycle).  

For example, the cycle $\pgap(\pml{5})$ has length (number of gaps) $\phi(\pml{5})=8$ and span
(sum of gaps) $\pml{5}=30$.
$$ \pgap(\pml{5}) \biggap = \biggap 6 \; 4 \; 2 \; 4 \; 2 \; 4 \; 6 \; 2.$$

There is a recursion from one cycle $\pgap(\pml{p_k})$ to the next $\pgap(\pml{p_{k+1}})$.  We concatenate $p_{k+1}$ copies of
$\pgap(\pml{p_k})$ and then add adjacent gaps at the running sums given by the elementwise product $p_{k+1}\ast \pgap(\pml{p_k})$. 
These additions of adjacent gaps are called {\em fusions}.  The gaps that survive the fusions are the gaps between primes.

The recursion across the cycles of gaps $\pgap(\pml{p_k})$ is a discrete dynamic system.  If we take initial conditions from 
$\pgap(\pml{p_0})$, then we can create {\em exact} population models for all gaps $g < 2 p_1$.  The populations all grow 
superexponentially, so we divide by factors of $p_k-2$ to obtain the exact models for relative populations for all gaps $g < 2p_1$.

\begin{eqnarray*}
w_g(\pml{p_{k+1}}) \; = \; \left[ \begin{array}{c} 
w_{g,1} \\ w_{g,2} \\ \vdots \\ w_{g,J} 
\end{array} \right]_{\pml{p_{k+1}}}  & = & 
 \left[ \begin{array}{ccccc} 
 1 & \lil \frac{1}{p_{k+1}-2} & 0 & \cdots & 0 \\
 0 & \lil \frac{p_{k+1}-3}{p_{k+1}-2} & \lil \frac{2}{p_{k+1}-2} & & 0 \\
 0 & 0 & \lil \frac{p_{k+1}-4}{p_{k+1}-2} & \ddots & 0 \\
 \vdots & \vdots & \ddots & \ddots & \vdots \\
 0 & 0 & \cdots & \cdots & \lil \frac{p_{k+1}-J-1}{p_{k+1}-2}
 \end{array} \right]
\cdot \left[ \begin{array}{c}
w_{g,1} \\ w_{g,2} \\  \vdots \\ w_{g,J} 
 \end{array} \right]_{\pml{p_k}}  \\[0.2cm]
 & = & M_J(p_{k+1}) \cdot w_g(\pml{p_k}) \\[0.2cm]
 & = & M_J^{k+1} \cdot w_g(\pml{p_0})
\end{eqnarray*}

Here the elements $w_{g,j}$ denote the relative population of driving terms for gap $g$ of length $j$.  $w_{g,1}(\pml{p})$ is the 
relative population of the gap $g$ itself in the cycle $\pgap(\pml{p})$.  And we use the notation $M^k$ to denote the product
of matrices:
$$
 M_J^k \; = \; M_J(p_k) \cdot M_J(p_{k-1}) \cdots M_J(p_2) \cdot M_J(p_1)
$$

$w_g(\pml{p_0})$ is the vector of the initial conditions in the cycle of gaps $\pgap(\pml{p_0})$.  These are the counts of the gap $g$
in this cycle and of all of its driving terms, divided by the population of gaps $2$.  
We need $J$ to be at least as large as the longest driving terms for $g$.

For more details about the discrete dynamic system and the population models, please see the prior work \cite{FBHSFU, FBHPatterns}.

\section{Models for $g < 2p_1$}
The iterative model
$$ w_g(\pml{p_k}) \; = \; M_J^k w_g(\pml{p_0})$$
only applies to gaps $g < 2 p_1$.

This constraint $g < 2p_1$ arises from needing to be sure that under the recursion each fusion occurs in its own copy of a driving 
term for $g$.  This allows us to get the exact counts across driving terms of all lengths.  
Since the fusions are spaced according to the elementwise product $p_{k+1} \ast \pgap(\pml{p_k})$ and the 
smallest element in $\pgap(\pml{p_k})$ is $2$, the fusions are separated by at least $2 p_k$.  So $g < 2p_0$ suffices to 
use the iterative system for the gap $g$.

The challenge to us developing models for larger gaps $g$ is that we need to get the initial populations from a cycle of gaps
$\pgap(\pml{p_0})$ such that $g < 2p_1$.  The cycle $\pgap(\pml{p})$ has length $\phi(\pml{p})$.

\begin{center}
\begin{tabular}{|c|cccccccc|} \hline
$p$ & $29$ & $31$ & $37$ & $41$ & $43$ & $47$ & $53$ & $59$ \\
$\max g$ & $60$ & $72$ & $80$ & $84$ & $92$ & $104$ & $116$ & $120$ \\
$\phi(\pml{p})$ & $\lil 1.02E9$ & $\lil 3..07E10$ & $\lil 1.10E12$ & $\lil 4.41E13$ & $\lil 1.85E15$ & $\lil 8.53E16$ & 
 $\lil 4.44E18$ & $\lil 2.57E20$ \\ \hline
\end{tabular}
\end{center}

We see that computing $\pgap(\pml{59})$ is at the horizon of current computing capability, and this would enable us
to calculate the initial conditions for the models for all gaps $g \le 120$.

For the models in our prior work \cite{FBHSFU, FBHPatterns} we used $\pgap(\pml{37})$ and exhibited the relative population
models for $g \le 80$.

\section{Models for $g=2p_1$}

We can extend the models to cover the case $g = 2p_1$.  The important insight here is that the populations of driving terms
for this gap $g$ can be exactly produced from the initial conditions $w_g(\pml{p_0})$.  We have to use a different system matrix
to update the relative populations for $g$ and its driving terms from $w_g(\pml{p_0})$ to $w_g(\pml{p_1})$.  After this one
special iteration we can apply the general model described above.

The special first iteration is motivated by the following Lemma.

\begin{lemma}
Let $s$ be a constellation of span $g$ in $\pgap(\pml{p_k})$.  If $p_{k+1} | g$, then both ends of $s$ are fused in the same
copy under the recursion from $\pgap(\pml{p_k})$ to $\pgap(\pml{p_{k+1}})$.
\end{lemma}

This is Lemma~28 on page~139 of \cite{FBHPatterns}

Let $s$ be a driving term for $g$ of length $j$. Under the recursion from $\pgap(\pml{p_0})$ to $\pgap(\pml{p_1})$, 
the concatenation step initially produces $p_1$ copies of $s$.  Each of the possible $j+1$ fusions in $s$ will occur
exactly once.  The $j-1$ interior fusions will result in a shorter driving term for $g$ in $\pgap(\pml{p_1})$, and the $2$
boundary fusions will eliminate that image of $s$ as a driving term for $g$.

In order to get a count of the driving terms for $g=2p_1$ in $\pgap(\pml{p_1})$, we need to confirm that the $j-2$ interior
fusions occur in separate images of $s$, and that the $2$ boundary fusions eliminate only $1$ image of $s$.

For any interior fusion in $s$, the span from this fusion to either end of $s$ is strictly less than $2p_1$.  But the smallest
distance between spans is $2p_1$, and thus the interior fusions will occur in separate images of $s$.  The interior fusions
in $s$ result in $j-1$ driving terms for $g$ of length $j-1$.

By the lemma above, the two boundary fusions occur in the same image of $s$, eliminating exactly one image of $s$
as a driving term for $g$.  Note that this is specific to $g=2p_1$; for $g < 2p_1$ the two boundary fusions occur in separate images
of $s$.

So for the gap $g=2p_1$, each driving term $s$ of length $j$ in $\pgap(\pml{p_0})$ produces $j-1$ driving terms of length
$j-1$ in $\pgap(\pml{p_1})$, one image of $s$ is removed as a driving term for $g$, and $p_1-j$ images of $s$ are preserved
intact as driving terms for $g$ in $\pgap(\pml{p_1})$.  This specific iteration for $w_{g=2p_1}$ is
\begin{eqnarray*}
w_g(\pml{p_1}) & = & \widehat{ M}_J(p_1) \cdot w_g(\pml{p_0}) \\[0.2cm]
 & = & 
 \left[ \begin{array}{cccccc} 
 \lil \frac{p_1-1}{p_1-2} & \lil \frac{1}{p_1-2} & 0 & 0 & \cdots & 0 \\
 0 & 1 & \lil \frac{2}{p_1-2} & 0 & & 0 \\
 0 & 0 & \lil \frac{p_1-3}{p_1-2} & \lil \frac{3}{p_1-2} & \ddots & 0 \\
 \vdots & \vdots & \ddots & \ddots & \ddots &  \vdots \\
 0 & 0 & 0 & 0 & \cdots & \lil \frac{p_1-J}{p_1-2}
 \end{array} \right] \cdot w_g(\pml{p_0})
\end{eqnarray*}

If we have computed the cycle $\pgap(\pml{p_0})$ for initial conditions, the counts of the driving terms across a range of gaps $g$,
we normalize these counts by the number of gaps $2$ and from these relative populations $w_g(\pml{p_0})$ we have the complete 
exact models for the relative populations across all subsequent stages of Eratosthenes sieve:
$$
w_g(\pml{p_k}) \; = \; \left\{ \begin{array}{cc}
M_J(p_k) \cdot M_J(p_{k-1}) \cdots M_J(p_2) \cdot M_J(p_1) \cdot w_g(\pml{p_0}) & {\rm if} \; g < 2p_1 \\[0.3cm]
M_J(p_k) \cdot M_J(p_{k-1}) \cdots M_J(p_2) \cdot \widehat{M}_J(p_1) \cdot w_g(\pml{p_0}) & {\rm if} \; g = 2p_1 
\end{array} \right.
$$

Specifically, from our initial conditions from $\pgap(\pml{37})$ we have previously been able to exhibit the population models for
all gaps $g \le 80$.  Here $p_0=37$ and $p_1=41$.  We can now add one more model, the model for $g=82$.

To harmonize the collection of models $w_g$ for $g \le 82$, we need the same starting point $p_0$.  We could use the dynamic system
to advance the models for all of the $g \le 80$ up to $w_g(\pml{41})$ and use $p_0=41$ as the starting point; or we could back
$w_{82}(\pml{41})$ up using $M^{-1}_J(41)$ to obtain an equivalent surrogate starting point $\widehat{w}_{82}(\pml{37})$ that
provides the exact model $w_{82}(\pml{p})$ for all $p \ge 41$.

We pursue this second approach here.
\begin{eqnarray*}
\widehat{w}_{82}(\pml{37}) & = & M^{-1}_J(41) \cdot w_{82}(\pml{41}) \\
  & = & M^{-1}_J(41) \cdot  \widehat{M}_J(41) \cdot w_{82}(\pml{37})
\end{eqnarray*}
This gives us a starting point $\widehat{w}_{82}(\pml{37})$ that aligns with our other starting points in $\pgap(\pml{37})$ and that
provides the exact values for $w_{82}(\pml{p_k})$ for all $p_k \ge 41$.  This surrogate starting point $\widehat{w}_{82}(\pml{37})$
is {\em not} the correct relative population $w_{82}(\pml{37})$.
\begin{eqnarray*}
\widehat{w}_{82}(\pml{37}) & \neq & w_{82}(\pml{37}) \\
M_J(41) \cdot \widehat{w}_{82}(\pml{37}) & = & w_{82}(\pml{41}) \; = \; \widehat{M}_J(41) \cdot w_{82}(\pml{37}) 
\end{eqnarray*}

Using $\pgap(\pml{37})$ for the starting point, we have driving terms for $g=82$ up to length $J=19$.  The number of gaps
$g=2$ in $\pgap(\pml{37})$ is 
$$n_{2,1}(\pml{37}) =217929355875.$$
For the gap $g=82$ we tabulate the data from $\pgap(\pml{37})$.  Our calculations have some numerical errors 
on the order of $10^{-13}$.

\begin{center}
\begin{table}
\begin{tabular}{|r|rrrrr|} \hline
$j$ & $n_{82,j}(\pml{37})$ & $w_{82,j}(\pml{37})$ & $w_{82,j}(\pml{41})$ & $\widehat{w}_{82,j}(\pml{37})$ & $l_{82,j}$ \\ \hline
$1$ &  $\lil 0$ & $\lil 0$ & $\lil 0$ & $\lil -2.768491E-13$  & $\lil 1.025641$ \\ 
$2$ &  $\lil 0$ & $\lil 0$ & $\lil 2.353150E-13$ & $\lil 1.976531E-12$  & $\lil 11.553942$ \\ 
$3$ &  $\lil 1$ & $\lil 4.5886429E-12$ & $\lil 1.160809E-09$ & $\lil -3.035074E-12$  & $\lil 60.410483$ \\ 
$4$ &  $\lil 3276$ & $\lil 1.5032394E-08$ & $\lil 1.415302E-07$ & $\lil 1.503443E-08	$ & $\lil 194.488465$   \\ 
$5$ &  $\lil 270422$ & $\lil 1.2408700E-06$ & $\lil 5.882215E-06$ & $\lil 1.244637E-06$ & $\lil 431.258857$  \\ 
$6$ &  $\lil 8051838$ & $\lil 3.6947010E-05$ & $\lil 1.179125E-04$ & $\lil 3.716878E-05$ & $\lil 697.937481$  \\ 
$7$ &  $\lil 120058788$ & $\lil 5.5090691E-04$ & $\lil 1.326320E-03$ & $\lil 5.558081E-04$ & $\lil 852.214506$ \\ 
$8$ &  $\lil 1027245782$ & $\lil 4.7136641E-03$ & $\lil 9.081749E-03$ & $\lil 4.769260E-03$ & $\lil 800.382107$  \\ 
$9$ &  $\lil 5411112020$ & $\lil 2.4829661E-02$ & $\lil 3.968207E-02$ & $\lil 2.519649E-02$ & $\lil 584.027075$  \\ 
$10$ & $\lil 18234669494$ & $\lil 8.3672387E-02$ & $\lil 1.136091E-01$ & $\lil 8.516773E-02$ & $\lil 332.122850$  \\ 
$11$ & $\lil 40031677310$ & $\lil 1.8369107E-01$ & $\lil 2.155096E-01$ & $\lil 1.875723E-01$ & $\lil 146.758457$ \\ 
$12$ & $\lil 57338080360$ & $\lil 2.6310398E-01$ & $\lil 2.702203E-01$ & $\lil 2.695709E-01$ & $\lil 49.939794$ \\ 
$13$ & $\lil 52822037198$ & $\lil 2.4238147E-01$ & $\lil 2.204693E-01$ & $\lil 2.492174E-01$ & $\lil 12.883935$  \\ 
$14$ & $\lil 30369623454$ & $\lil 1.3935536E-01$ & $\lil 1.135898E-01$ & $\lil 1.438024E-01$ & $\lil 2.461383$ \\ 
$15$ & $\lil 10389093440$ & $\lil 4.7671840E-02$ & $\lil 3.526600E-02$ & $\lil 4.936702E-02$ & $\lil 0.336767$ \\ 
$16$ & $\lil 1974527214$ & $\lil 9.0604004E-03$ & $\lil 6.171214E-03$ & $\lil 9.413221E-03$ & $\lil 0.031541$ \\ 
$17$ & $\lil 192967582$ & $\lil 8.8545933E-04$ & $\lil 5.642306E-04$ & $\lil 9.225035E-04$ & $\lil 0.001910$ \\ 
$18$ & $\lil 9665424$ & $\lil 4.4351180E-05$ & $\lil 2.673245E-05$  & $\lil 4.631847E-05$ & $\lil 0.000070$ \\ 
$19$ & $\lil 272272$ & $\lil 1.2493590E-06$ & $\lil 7.047666E-07$ & $\lil 1.308852E-06$ & $\lil 0.000001$ \\ \hline
  sum & $\lil 217929355875$ & $\lil 1.000000$ & $\lil 1.025641$ & $\lil 1.025641$ & \\ \hline
\end{tabular}
\caption{\label{w82Tbl} Data for the gap $g=82$ from the cycle $\pgap(\pml{37})$.  The column $n_{82,j}$ lists the
populations of driving terms for $g=82$ for various lengths $j$.  The longest driving term has length $J=19$.  The column 
$w_{82,j}(\pml{37})$ contains the normalized populations for these driving terms, and $w_{82,j}(\pml{41})$ calculates
the relative populations for the driving terms in $\pgap(\pml{41})$ using $\widehat{M}_{19}(41)$.  The next column
$\widehat{w}_{82,j}(\pml{37})$ is the pre-image of $w_{82,j}(\pml{41})$ under $M^{-1}_{19}(41)$.  The final column
lists the coefficients $l_{82,j}$ for the population model. }
\end{table}
\end{center}

The models update the relative populations of all of the driving terms for $g=82$.  If we take just the top row we extract the
model just for the relative population of the gap $g=82$ itself.
$$
w_{82,1}(\pml{p_k})  = \; l_{82,1} - l_{82,2} \prod_{41}^{p_k} \frac{q-3}{q-2} +  l_{82,3} \prod_{41}^{p_k} \frac{q-4}{q-2} 
 - \cdots +  l_{82,19} \prod_{41}^{p_k} \frac{q-20}{q-2} 
$$
for all $p_k \ge 41$.
Table~\ref{w82Tbl} lists these coefficients $l_{82,j}$ in the final column.

\section{Notes on models for slightly larger $g$}

The beauty of the work above is that we can work directly with the relative populations of driving terms of various
lengths $j$.  Can we extend these methods any further, to extract models for $g=2p_1+2$ or $g=2p_1+4$?
To do so, we would need to track subpopulations among the driving terms.

\begin{figure}[htb]
\centering
\includegraphics[width=5in]{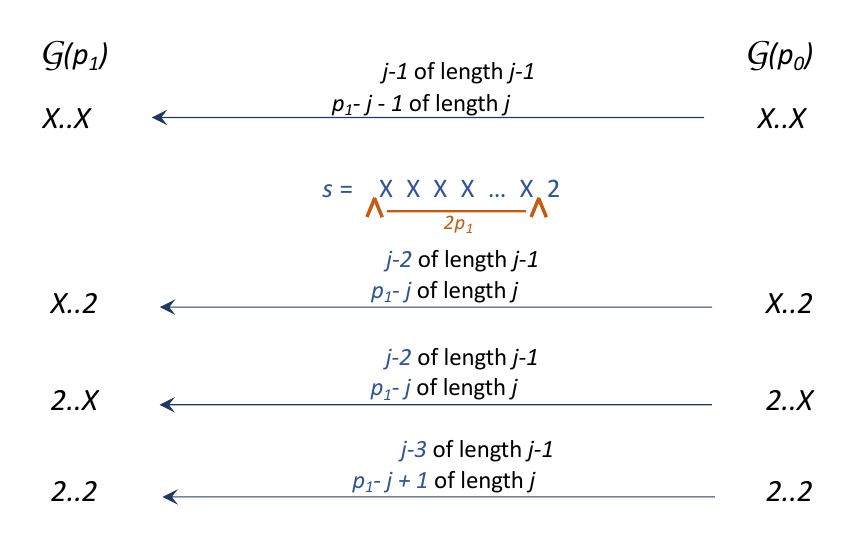}
\caption{\label{2p2Fig} A diagram of the specific iteration for gaps $g=2p_1+2$ from $\pgap(\pml{p_0})$ to $\pgap(\pml{p_1})$.
'X' denotes any gap in a driving term $s$ for $g$ that is not a $2$.  If the driving term $s$ does not begin or end with a $2$, then
the general model applies.  If $s$ ends with a $2$, then that interior fusion occurs in the same image of $s$ as the far boundary
fusion.  One fewer image of length $j-1$ survives, and one more image of length $j$ survives.}
\end{figure}

{\em Gaps $g=2p_1+2$.}  This gap is small enough that $g < 2p_2$, so the only spans that complicate our tracking
the counts of the driving terms for $g$ are the spans of $2p_1$ where the two fusions occur as a boundary fusion and
interior fusion in a single image of a driving term $s$.  This will occur iff the driving term $s$ begins or ends with a $2$.

If we separate the counts for the driving terms for $g=2p_1+2$ to cover four subpopulations, we could perform a distinct
iteration from $w_{2p_1+2}(\pml{p_0})$ to $w_{2p_1+2}(\pml{p_1})$, after which we can use the general model. 

Consider a driving term $s$ for $g$ of length $j$ in $\pgap(\pml{p_0})$.  Since $g=2p_1+2$ the two boundary fusions
occur in different images of $s$, so $p_1-2$ images of $s$ survive as driving terms for $g$ in $\pgap(\pml{p_1})$.
We need to track their lengths.  If the first or last gap in $s$ is a $2$, then by the lemma this interior fusion occurs
in the same image of $s$ as the far boundary fusion. 

We have the following four subpopulations of driving terms $s$ of span $2p_1+2$ and length $j$ in $\pgap(\pml{p_0})$:
\begin{enumerate}[label=(\alph*)]
\item $s = X \ldots X$.  If $s$ begins and ends with gaps other than $2$, then all of the fusions - interior and boundary - occur
in separate images of $s$ during the recursive construction of $\pgap(\pml{p_1})$.  We can use the general model for these
populations described above.
\item $s = X \ldots 2$.  If $s$ ends with a $2$, then the interior fusion for this last gap will occur in the same image of $s$ as the
boundary fusion at the start of $s$.  Of the $p_1$ images of $s$ created during the concatenation step, two are eliminated as
driving terms for $g$ by the boundary fusions.  The one boundary fusion takes an interior fusion along with it.  The remaining
$j-2$ interior fusions result in driving terms for $g$ in $\pgap(\pml{p_1})$ of length $j-1$.  $p_1-j$ images of $s$ survive intact
as driving terms of length $j$ for $g$ in $\pgap(\pml{p_1})$.
\item $s = 2 \ldots X$.  The cycle of gaps $\pgap(\pml{p})$ is symmetric.  So if $s$ starts with a $2$, the analysis is the same
as the previous case.
\item $s = 2 \ldots 2$.  For driving terms $s$ that begin {\em and} end with a $2$, these would fall into both of the previous 
cases, complicating the counts.  So we separate them out as their own subpopulation.  In this case {\em both} boundary
fusions coincide with an interior fusion.  Of the $j-1$ interior fusions, $j-3$ produce driving terms for $g$ in $\pgap(\pml{p_1})$
of length $j-1$.   $p_1-j+1$ images of $s$ survive intact as driving terms of length $j$ for $g$ in $\pgap(\pml{p_1})$.
\end{enumerate}

We could exhibit the block banded matrices for this system.  But unless we partition the population of driving terms for the
gap $g=2p_1+2$ into the required subpopulations, we cannot apply this model.

{\em Gaps $g=2p_1+4$.} 
For gaps $g=2p_1+4$, the subpopulations of driving terms become more complicated.  We need to consider
cases in which a driving term begins or ends with a $4$, and the analysis parallels the work above for $g=2p_1+2$.

We also have to consider the possibility that $g=2p_2$, which would occur when $p_2=p_1+2$.

\section{Conclusion}
This work serves as an addendum to the existing references \cite{FBHSFU,FBHPatterns}.  
We do not duplicate that background here, beyond summarizing a few needed results.

We have shown previously that at each stage of Eratosthenes sieve there is a corresponding cycle of gaps $\pgap(\pml{p})$.
We can view these cycles of gaps as a discrete dynamic system, and from this system we can obtain exact models for 
the populations and relative populations of gaps $g < 2p_1$ if we can get the initial conditions from $\pgap(\pml{p_0})$.

In this addendum we have shown that we can produce the model for ${g=2p_1}$ from these initial conditions.  This model requires
one special iteration to track the count from $\pgap(\pml{p_0})$ to $\pgap(\pml{p_1})$, after which we can use the general model for these populations.

We have further shown that in order to produce the models for ${g=2p_1+2}$ and beyond from initial conditions in $\pgap(\pml{p_0})$,
we would have to track subpopulations of the driving terms $s$ until the general model applies, that is until ${g < 2p_{k+1}}$.


\bibliographystyle{amsplain}

\providecommand{\bysame}{\leavevmode\hbox to3em{\hrulefill}\thinspace}
\providecommand{\MR}{\relax\ifhmode\unskip\space\fi MR }
\providecommand{\MRhref}[2]{%
  \href{http://www.ams.org/mathscinet-getitem?mr=#1}{#2}
}
\providecommand{\href}[2]{#2}

\end{document}